\documentclass[12pt]{amsart}

\usepackage{amsthm}
\usepackage{amsmath}
\usepackage{amssymb}
\usepackage{latexsym}
\usepackage{enumerate}

\usepackage{graphicx}

\renewcommand{\thetheoremName}

\newtheorem{theorem}{Theorem}[section]
\newtheorem{lemma}[theorem]{Lemma}

\newtheorem{proposition}[theorem]{Proposition}
\newtheorem{corollary}[theorem]{Corollary}

\theoremstyle{definition}
\newtheorem{definition}[theorem]{Definition}

\newtheorem{remark}[theorem]{Remark}

\numberwithin{equation}{section}

\newcommand{\Hess}{\operatorname{Hess}}
\newcommand{\dist}{\operatorname{dist}}
\newcommand{\Vol}{\operatorname{Vol}}

\newcommand{\Div}{\operatorname{div}}
\newcommand{\C}{\operatorname{Cap}}

\newcommand{\erre}{\mathbb{R}}
\newcommand{\Lmod}{\operatorname{L}}

\newcommand{\pL}{{\operatorname{\Delta}}_{p}}

\newcommand{\pLP}{{\operatorname{\Delta}}_{p}^{S}}

\newcommand{\pCap}{{\operatorname{Cap}}_{p}}
\newcommand{\loc}{{\rm loc}}
\begin{document}

\title[$p$-Parabolicity]{A note on the $p$-Parabolicity\\
of Submanifolds}

\author[A. Hurtado]{Ana Hurtado*}
\address{Departamento de Geometr\'\i a y Topolog\'\i a, Universidad de Granada, Granada,
Spain.} \email{ahurtado@ugr.es}
\thanks{* Work partially supported by
the Caixa Castell\'{o} Foundation, DGI grant MTM2007-62344 and
Junta de Andaluc\'{\i}a grants P06-FQM-01642 and FQM325.}

\author[V. Palmer]{Vicente Palmer*}
\address{Departament de Matem\`{a}tiques-Institut de Matem\'atiques i Aplicacions de Castell\'o, Universitat Jaume I, Castellon,
Spain.} \email{palmer@mat.uji.es}
\thanks{* Work partially supported by
the Caixa Castell\'{o} Foundation, and  DGI grant MTM2007-62344.}

%    General info
%\date{January 1, 1994 and, in revised form, June 22, 1994.}

%\dedicatory{This paper is dedicated to [[[]]].}


\maketitle
\bibliographystyle{acm}
\begin{abstract}
We give a geometric criterion which shows $p$-parabolicity of a class of submanifolds in a Riemannian manifold, with controlled second fundamental form, for $p \geq 2$.
\end{abstract}

%%%%%%%%%%%%%%%%%%%%%%%%%%%%%%%%%%%%%%%%%%%%%%%%%%%%%%%%%%%%%%%%%%%%%%%%
%
%       SECTION 1 Introduction
%
%%%%%%%%%%%%%%%%%%%%%%%%%%%%%%%%%%%%%%%%%%%%%%%%%%%%%%%%%%%%%%%%%%%%%%%%

\section{Introduction} \label{secIntro}
A Riemannian manifold $M^n$ is called {\em parabolic} if it does not admit a nonconstant positive superharmonic function (referred to the $2$-Laplacian). Otherwise is called {\em hyperbolic}. In the paper \cite{LS}, T. Lyons and D. Sullivan stablished a list of equivalent conditions to check the  $2$-hyperbolicity, (and hence, the $2$- parabolicity) of an oriented Riemannian manifold, under the name of {\em the Kelvin-Nevanlinna-Royden criterion}. Following this research, in the paper \cite{GT}, V. Gol'dshtein and M. Troyanov extended this criterion to the $p$-Laplace operator, ($1 <p< \infty$).

Among the characterizations given by the Kelvin-Nevanlinna-Royden
criterion for $p$-parabolicity of a Riemannian manifold, we shall
use in this paper those that states that a manifold $M$ is
$p$-parabolic, (with $1 <p< \infty$), if and only if it has
vanishing $p$-capacity, namely, there exists a compact set $D
\subseteq M$ such that $Cap_p(D,M)=0$.

Here, the $p$-capacity of $D$ is defined by
$$Cap_p(D,M)=\inf_{u}\int_M \Vert\nabla u\Vert^p d\mu,$$
where the infimum is taken over all real-valued functions $u\in C^{\infty}_0(M)$, with $u\geq 1$ in $D$, (see \cite{HKM}).

On the other hand, in \cite{GT} and \cite{T1}  it is possible to
find a set of geometric criteria for this so-called {\em type
problem} for the $p$-Laplacian in Riemannian manifolds, namely, to
decide when a Riemannian manifold is $p$-parabolic of
$p$-hyperbolic.

In particular, in the Corollary 5.2 of \cite{T1}, (see too Proposition 4 in \cite{GT}), it is presented a characterization of $p$-parabolicity for Riemannian manifolds with a warped cilindrical end, and in Corollary 5.4 of \cite{T1}, (see too Proposition 3 in \cite{GT}), we have a sufficient codition for $p$-parabolcity in terms of the volume growth of the manifold.

While these two criteria are intrinsic, we are going to present in this paper a geometric criterion to decide if a submanifold $S^m$ properly immersed in an ambient manifold $N^n$ with a pole is $p$-parabolic, which involves (lower)  bounds for the mean curvature and the second fundamental form of $S$.

This criterion is based, (as in \cite{MP1} and, specially, \cite{HMP} from which this paper can be considered a spin-off), in the Hessian-Index analysis of the (restricted to $S$) extrinsic distance function from the pole, (see \cite{GreW}).

%%%%%%%%%%%%%%%%%%%%%%%%%%%%%%%%%%%%%%%%%%%%%%%%%%%%%%%%%%%%%%%%%%%%%%%%
%
%       SUBSECTION  Outline
%
%%%%%%%%%%%%%%%%%%%%%%%%%%%%%%%%%%%%%%%%%%%%%%%%%%%%%%%%%%%%%%%%%%%%%%%%
\subsection{Outline of the paper} \label{subsecOutline}
We shall present the basic definitions concerning the
$p$-Laplacian in Section 2. Section 3 is devoted to the study of
the curvature setting where our results hold, together with the
Hessian and Laplacian analysis needed. Main results are stated and
proved in Sections 4, 5 and 6.

\subsection{Acknowledgements}
We would like to acknowledge professors Ilkka Holopainen and Steen
Markvorsen their useful comments concerning these results.

%%%%%%%%%%%%%%%%%%%%%%%%%%%%%%%%%%%%%%%%%%%%%%%%%%%%%%%%%%%%%%%%%%%%%%%%
%
%       SECTION 2 p-Laplacian
%
%%%%%%%%%%%%%%%%%%%%%%%%%%%%%%%%%%%%%%%%%%%%%%%%%%%%%%%%%%%%%%%%%%%%%%%%

\section{The $p$-Laplacian}\label{secpLap}
Let $M$ be a non-compact Riemannian manifold, with the Riemannian metric
$\langle\cdot,\cdot\rangle$ and the Riemannian volume form $d\mu$.

The $p$-Laplacian of a $C^2$ function $u$ is defined as
$$\Delta_p u = \Div(\Vert \nabla u\Vert^{p-2} \nabla u).$$

When $p=2$, we have the usual Laplacian, and the classical potential theory developed from the study of the solutions of the Laplace equation
$$\Delta_2 u=0.$$
However, when $p \neq 2$, equation
\begin{equation}\label{pLaplaceEq}
\Delta_pu=0,
\end{equation}
is nonlinear and degenerates at the zeroes of the gradient of $u$.

Then, the solutions of (\ref{pLaplaceEq}), known usually as
{\em $p$-harmonic functions}, need not be smooth, nor even $C^2$ and
equation (\ref{pLaplaceEq}) must be interpreted in a {\em weak} sense.

In this way, and given $1<p<\infty$, we say that a function $u \in W^{1,p}_{\loc}(M)$ is a (weak) solution to the
$p$-Laplace equation
\begin{equation} \label{eqPlap}
-\Div\left( \Vert \nabla u\Vert^{p-2} \nabla u
\right)  =   0
\end{equation}
in $M$ if
\begin{equation} \label{eqPlapInt}
\int_{M} \langle\Vert \nabla u\Vert^{p-2} \nabla u, \nabla
\phi \rangle \, d\mu  =   0
\end{equation}
for all $\phi \in C^{\infty}_{0}(M)$.

Here, $\nabla u\in L^{1}_{\loc}(M)$ is the {\em{distributional gradient}} of
 $u\in L^{1}_{\loc}(M)$. Furthermore, $L^1(M)$ denotes the space of measurable functions $f:M \longrightarrow \erre$ with finite norm $\Vert f\Vert_1 <\infty$, and $L^{1}_{\loc}(M)$ is its corresponding local space defined through the open sets in $M$ with compact closure, (see \cite{HKM} p. 13). In its turn, the space
$W^{1,p}(M),\ 1\le p<\infty$ is the Sobolev space of all functions
$u\in L^{p}(M)$ whose distributional gradient $\nabla u$ belongs to
$L^{p}(M)$, equipped with  the norm
$\Vert u\Vert_{1,p}=\Vert u\Vert_{p}+\Vert\nabla u\Vert_{p}$.
The corresponding local space $W^{1,p}_{\loc}(M)$ is defined as in \cite{HKM}.

Continuous
solutions of \eqref{eqPlap} are called
{\em{$p$-harmonic}}. Here the continuity
assumption makes no restriction since every
solution of \eqref{eqPlap} has a continuous
representative. The extension of regularity results of this kind, (see \cite{E} and \cite{Li}), from the Euclidean setting to the Riemannian setting is detailed in \cite{HMP}, Remark 9.2.

A function $u \in W^{1, p}_{\loc}(M)$ is called a
{\em{$p$-supersolution}} in $M$ if
\begin{equation*}
\int_{M} \langle \Vert \nabla u\Vert^{p-2} \nabla u, \nabla \phi
\rangle \, d\mu  \geq  0,
\end{equation*}
for all non-negative $\phi \in C^{\infty}_{0}(M)$. If, moreover,
$u$ is lower semicontinuous, then $u$ is $p$-superharmonic, (and
we denote $\Delta_p u \leq 0$).

Similarly, a function $v \in W^{1, p}_{\loc}(M)$ is called a
{\em{$p$-subsolution}} in $M$ if $-v$ is a $p$-supersolution.
 If, moreover, $v$ is lower semicontinuous,
then $v$ is $p$-subharmonic, ($\Delta_p v \geq 0$).

A fundamental feature of
solutions of \eqref{eqPlap} is the following well-known {\em{maximum (or comparison)
principle}} which will be instrumental for the comparison technique
presented below in Sections \ref{proofthmMain1} and \ref{proofthmMain2}:
If $u\in W^{1,p}(M)$ is a $p$-supersolution,
$v\in W^{1,p}(M)$ is a $p$-subsolution, and $\max(v-u,0)\in W^{1,p}_{0}(M)$,
then $u\ge v$ a.e. in $M$. In particular, if $D\subset M$ is a
precompact open set, $u\in C(\bar D)$ is a $p$-supersolution,
$v\in C(\bar D)$ is a $p$-subsolution, and $u\ge v$ in $\partial D$,
then $u\ge v$ in $D$.
We refer to \cite[3.18]{HKM} for a short proof of the comparison
principle.
%For the reader's convenience we recall the short proof
%of the comparison principle from \cite[3.18]{HKM}. The proof is based on the
%following elementary inequality:
%Let $a \neq b$ denote two vectors in a given tangent space $T_{x}M$
%and suppose that $1 < p < \infty$. Then
%\[
%\langle \Vert a \Vert^{p-2}a - \Vert b \Vert^{p-2}b,  a - b
%\rangle >  0.
%\]
%Suppose then that $u\in W^{1,p}(M)$ is a $p$-supersolution and
%$v\in W^{1,p}(M)$ is a $p$-subsolution such that
%$\phi=\max(v-u,0)\in W^{1,p}_{0}(M)$.
%Since
%\begin{align*}
%0&\ge \int_{M}\langle\Vert\nabla v\Vert^{p-2}\nabla
%v,\nabla\phi\rangle\,d\mu
%-\int_{M}\langle\Vert\nabla u\Vert^{p-2}\nabla
%u,\nabla\phi\rangle\,d\mu\\
%&=\int_{\{u<v\}}\langle\Vert\nabla v\Vert^{p-2}\nabla v
%-\Vert\nabla u\Vert^{p-2}\nabla u,\nabla v-\nabla u\rangle\,d\mu\ge 0,
%\end{align*}
%we have $\nabla\phi =0$ a.e. in $M$ as required.

%%%%%%%%%%%%%%%%%%%%%%%%%%%%%
%   Section 3: Comparison Constellations
%%%%%%%%%%%%%%%%%%%%%%%%%

\section{Comparison Constellations} \label{secCompConstel}
We assume throughout the paper that $S^{m}$ is a non-compact,
properly immersed, and connected  Riemannian submanifold of a
complete Riemannian manifold $N^{n}$. Furthermore, we assume that
$N^{n}$ possesses at least one pole. Recall that a pole is a point
$o$ such that the exponential map $\exp_{o}\colon T_{o}N^{n} \to
N^{n}$ is a diffeomorphism. For every $x \in N^{n}\setminus \{o\}$
we define $r(x) = \dist_{N}(o, x)$, and this distance is realized
by the length of a unique geodesic from $o$ to $x$, which is the
{\it radial geodesic from $o$}. We also denote by $r$ the
restriction $r\vert_S: S\to \erre_{+} \cup \{0\}$. This
restriction is called the {\em{extrinsic distance function}} from
$o$ in $S^m$. The gradients of $r$ in $N$ and $S$ are denoted by
$\nabla^N r$ and $\nabla^S r$, respectively. Let us remark that
$\nabla^S r(x)$ is just the tangential component in $S$ of
$\nabla^N r(x)$, for all $x\in S$. Then we have the following
basic relation:
\begin{equation*}
\nabla^N r = \nabla^S r +(\nabla^N r)^\bot ,
\end{equation*}
where $(\nabla^N r)^\bot(x)$ is perpendicular to
$T_{x}S$ for all $x\in S$.

With the extrinsic distance at hand, we define the following domains:

\begin{definition}
Given a connected and complete
$m$-dimen\-sional submanifold $S^m$ in a complete
Riemannian manifold $N^n$ with a pole $o$, we
denote
the {\em{extrinsic metric balls}} of
(sufficiently large) radius $R$ and center $o$ by
$D_R(o)$. They are defined as any connected
component of the intersection
$$
B_{R}(o) \cap S =\{x\in S \colon r(x)< R\},
$$
where $B_{R}(o)$ denotes the open geodesic ball
of radius $R$ centered at the pole $o$ in
$N^{n}$.
Using these extrinsic balls we define the
$o$-centered extrinsic annuli
$$
A_{\rho,R}(o)= D_R(o) \setminus \bar D_{\rho}(o)
$$
in $S^m$ for $\rho < R$, where $D_{R}(o)$ is the component
of $B_{R}(o) \cap S$ containing $D_{\rho}(o)$.
\end{definition}

\begin{remark}\label{theRemk0}
We must to point out that these extrinsic domains are precompact,
(because the submanifold $S$ is properly immersed), and that the
radii $R$ that produce smooth boundaries $\partial D_{R}(o)$ are
dense in $\mathbb{R}$ by Sard's theorem and the Regular Level Set
Theorem, because the distance function $r$ is smooth in
$N^{n}\setminus \{o\}$, and hence, its restriction to $S$,
$r\vert_S$.
\end{remark}
%%%%%%%%%%%%%%%%%%%%%%%%%%%%%%%
%  Subsection: Curvature restrictions
%%%%%%%%%%%%%%%%%%%%%%%%%%

\subsection{Curvature restrictions}
\label{subsecCurvRestrict}

\begin{definition}\label{defRadCurv}
Let $o$ be a point in a Riemannian manifold $M$
and let $x \in M\setminus\{ o \}$. The sectional
curvature $K_{M}(\sigma_{x})$ of the two-plane
$\sigma_{x} \in T_{x}M$ is then called an
\textit{$o$-radial sectional curvature} of $M$ at
$x$ if $\sigma_{x}$ contains the tangent vector
to a minimal geodesic from $o$ to $x$. We denote
these curvatures by $K_{o, M}(\sigma_{x})$.
\end{definition}

\begin{definition}
The {\em{$o$-radial mean convexity}}
$\,\mathcal{C}(x)$ of $S$ in $N$, is defined  as
follows:
\begin{equation*}
\mathcal{C}(x) \,=\, -\langle \nabla^{N}r(x),
H_{S}(x) \rangle, \quad x \in S,
\end{equation*}
where $H_{S}(x)$ denotes the mean curvature
vector of $S$ in $N$.
\end{definition}

Moreover, for $p > 2$ we define:

\begin{definition}\label{defBetax}
The {\em{$o$-radial component $\,\mathcal{B}(x)$
of the second fundamental form}}  of $S$ in $N$,
is defined as:
\begin{equation*}
\mathcal{B}(x)  =  - \langle\nabla^{N}r(x),
\alpha_{x}(U_{r}, U_{r})\rangle,
\end{equation*}
where
\begin{equation*}
U_{r} =  \nabla^{S}(r(x))/\Vert
\nabla^{S}r(x) \Vert \in T_{x}S \subset T_{x}N
\end{equation*}
is the unit tangent vector to $S$ in the
direction of $\nabla^{S}r(x)$ (resp. tacitly
assumed to be $0$ in case $\nabla^{S}r(x) = 0$).
\end{definition}

Finally,

\begin{definition}
The {\em{\,$o$-radial tangency}} $\mathcal{T}(x)$
of $S$ in $N$ is defined as follows:
\begin{equation*}
\mathcal{T}(x)  =  \Vert \nabla^{S}r(x)\Vert
\end{equation*}
for all $x\in S$.
\end{definition}

Upper and lower bounds of $\mathcal{C}(x)$, $\mathcal{B}(x)$ and
$\mathcal{T}(x)$  together with a suitable control on the
$o$-radial sectional curvatures of the ambient space will
eventually control the $p$-Laplacian of restricted radial
functions on $S$.

%%%%%%%%%%%%%%%%%%
%Subsection: warped products and model spaces
%%%%%%%%%%%%%%%%%%%%

\subsection{Model spaces}\label{subsecWarp}

\begin{definition}[See \cite{Gri}, \cite{GreW}]\label{defModel}
A $w-$model $M_{w}^{m}$ is a smooth warped product with base $B^{1}
= [0,\Lambda[ \,\subset \mathbb{R}$ (where $0 < \Lambda
\leq  \infty$), fiber $F^{m-1} = \mathbb{S}^{m-1}_{1}$ (i.e. the unit
$(m-1)$-sphere with standard metric), and warping function $w\colon
[0,\Lambda[ \to \mathbb{R}_{+}\cup \{0\}$, with $w(0) = 0$,
$w'(0) = 1$, and $w(r) > 0$ for all $r >  0$. The point
$o_{w} = \pi^{-1}(0)$, where $\pi$ denotes the projection onto
$B^1$, is called the {\em{center point}} of the model space. If
$\Lambda = \infty$, then $o_{w}$ is a pole of $M_{w}^{m}$.
\end{definition}

\begin{remark}\label{propSpaceForm}
The simply connected space forms $\mathbb{K}^{m}(b)$ of constant
curvature $b$ are $w-$models as we mentioned in \cite{HMP}, (see \cite{GreW} and \cite{Gri}).
\end{remark}

\begin{proposition}[See \cite{GreW} and \cite{Gri}]\label{propWarpMean}
Let $M_{w}^{m}$ be a $w-$model with warping function $w(r)$ and
center $o_{w}$. The distance sphere of radius $r$ and center
$o_{w}$ in $M_{w}^{m}$ is the fiber $\pi^{-1}(r)$. This distance
sphere has the constant mean curvature $\eta_{w}(r)=
\frac{w'(r)}{w(r)}$. On the other hand, the $o_{w}$-radial
sectional curvatures of $M_{w}^{m}$ at every $x \in \pi^{-1}(r)$
(for $r > 0$) are all identical and determined by
\begin{equation*}
K_{o_{w} , M_{w}}(\sigma_{x}) = -\frac{w''(r)}{w(r)}.
\end{equation*}
\end{proposition}

%%%%%%%%%%%%%%%%%%
%Subsection: comparison constellations
%%%%%%%%%%%%%%%%%%%%

\subsection{Comparison constellations}\label{subsecCompConst}

We now collect the previous ingredients and
formulate the general framework for our
$p$-parabolicity comparison result, which results a dual setting
with respect the curvature assumptions stated in \cite{HMP} to obtain p-hyperbolicity.

\begin{definition}\label{defLowerConstellat}
Let $N^{n}$ denote a Riemannian manifold with a
pole $o$. Let $S^{m}$ denote a connected
complete submanifold properly immersed in $N^{n}$.
 Let $M_{w}^{m}$
denote a $w$-model with center
$o_{w}$; see Definition \ref{defModel}. We shall assume that
the
$o$-radial sectional curvatures of $N$ are
bounded from below by the $o_{w}$-radial
sectional curvatures of $M_{w}^{m}$:
\begin{equation} \label{eqKcomp}
K_{o, N}(\sigma_{x}) \geq
-\frac{w''(r)}{w(r)}
\end{equation}
for all $x$ with $r=r(x)\in [0,R]$.

Then the triple $\{ N^{n}, S^{m},
M_{w}^{m} \}$ is called a {\em{comparison
constellation with lower tangency}} on the interval $[0, R]$ if
the radial tangency $\mathcal{T}$ and the radial
convexity functions  $\mathcal{B}$ and
$\mathcal{C}$ of the submanifold $S^{m}$ are all
bounded from below by smooth radial functions
$g(r)$, $\lambda(r)$, and $h(r)$, respectively:
\begin{equation}
\label{eqRadConv}
\begin{aligned}
\mathcal{T}(x) &\geq g(r(x)),\\
\mathcal{B}(x) &\geq \lambda(r(x)) ,\,\,\, \textrm{and}\\
\mathcal{C}(x) &\geq h(r(x)) \,\,\,\textrm{for
all}\,\,\, x \in S^{m} \,\,\,
{\textrm{with}}\,\,\, r(x) \in [0, R] .
 \end{aligned}
\end{equation}
\end{definition}

Note that the radial tangency is, in a natural way, bounded from above by 1. This fact motivates the following
\begin{definition}
We assume the same general hypothesis on $S$ and $N$ than in the definition above.
The triple $\{ N^{n}, S^{m},
M_{w}^{m} \}$ is called a {\em{comparison
constellation with upper tangency}} on the interval $[0, R]$
when  the radial
convexity functions  $\mathcal{B}$ and
$\mathcal{C}$ of the submanifold $S^{m}$ are all
bounded from below by smooth radial functions
 $\lambda(r)$, and $h(r)$, respectively:
\begin{equation}
\label{eqRadConv}
\begin{aligned}
%\mathcal{T}(x) &\leq g(r(x)),\\
\mathcal{B}(x) &\geq \lambda(r(x)) ,\,\,\, \textrm{and}\\
\mathcal{C}(x) &\geq h(r(x)) \,\,\,\textrm{for
all}\,\,\, x \in S^{m} \,\,\,
{\textrm{with}}\,\,\, r(x) \in [0, R] .
 \end{aligned}
\end{equation}
\end{definition}

%%%%%%%%%%%%%%%%%%%%%%%%
%Subsection: Hessian and Laplacian comparison
%%%%%%%%%%%%%%%

\subsection{Hessian and Laplacian comparison analysis}\label{subsecLap}
The 2nd order analysis of the restricted distance function
$r_{|_{P}}$ defined on manifolds with a pole is firstly and
foremost governed by the Hessian comparison Theorem A in
\cite{GreW}:

\begin{theorem}[See \cite{GreW}, Theorem A]\label{thmGreW}
Let $N=N^{n}$ be a manifold with a pole $o$, let $M=M_{w}^{m}$
denote a $w-$model with center $o_{w}$, and $m \leq n$. Suppose
that every $o$-radial sectional curvature at $x \in N \setminus
\{o\}$ is bounded from below by the $o_{w}$-radial sectional
curvatures in $M_{w}^{m}$ as follows:
\begin{equation*}
K_{o, N}(\sigma_{x}) \geq -\frac{w''(r)}{w(r)}
\end{equation*}
for every radial two-plane $\sigma_{x} \in T_{x}N$ at distance $r =
r(x) = \dist_{N}(o, x)$ from $o$ in $N$. Then the Hessian of the
distance function in $N$ satisfies
\begin{equation}\label{eqHess}
\begin{aligned}
\Hess^{N}(r(x))(X, X) &\leq \Hess^{M}(r(y))(Y, Y)\\ &=
\eta_{w}(r)\left(1 - \langle \nabla^{M}r(y), Y \rangle_{M}^{2}
\right) \\ &= \eta_{w}(r)\left(1 - \langle \nabla^{N}r(x), X
\rangle_{N}^{2} \right)
\end{aligned}
\end{equation}
for every unit vector $X$ in $T_{x}N$ and for every unit vector $Y$
in $T_{y}M$ with $\,r(y) = r(x) = r\,$ and $\, \langle
\nabla^{M}r(y), Y \rangle_{M} = \langle \nabla^{N}r(x), X
\rangle_{N}\,$.
\end{theorem}

As a consequence of this result, we have the following Laplacian
inequality:
\begin{proposition} \label{corLapComp}
Let $N^{n}$ be a manifold with a pole $o$, let $M_{w}^{m}$ denote
a $w-$model with center $o_{w}$. Suppose that every $p$-radial
sectional curvature at $x \in N - \{o\}$ is bounded from below by
the $o_{w}$-radial sectional curvatures in $M_{w}^{m}$ as follows:
\begin{equation}\label{eqKbound}
\mathcal{K}(\sigma(x)) \, = \, K_{o, N}(\sigma_{x}) \geq
-\frac{w''(r)}{w(r)}
\end{equation}
for every radial two-plane $\sigma_{x} \in T_{x}N$ at distance $r
= r(x) = \dist_{N}(o, x)$ from $o$ in $N$. Then we have for every
smooth function $f(r)$ with $f'(r) \leq 0\,\,\textrm{for
all}\,\,\, r$, (respectively $f'(r) \geq 0\,\,\textrm{for
all}\,\,\, r$):
\begin{equation} \label{eqLap1}
\begin{aligned}
\Delta^{S}(f \circ r) \, \geq (\leq) \, &\left(\, f''(r) -
f'(r)\eta_{w}(r) \, \right)
 \Vert \nabla^{S} r \Vert^{2} \\ &+ mf'(r) \left(\, \eta_{w}(r) +
\langle \, \nabla^{N}r, \, H_{S}  \, \rangle  \, \right)  \quad ,
\end{aligned}
\end{equation}
where $H_{S}$ denotes the mean curvature vector of $S$ in $N$.
\end{proposition}

%%%%%%%%%%%%%%%%%%%%%%%%%%
%Section 4: Main results
%%%%%%%%%%%%%%%%%%%%%%

\section{Main results}\label{secMain}
Applying the notion of a comparison constellation
as defined in the previous section, we now
formulate our main $p$-parabolicity results. The
proofs are developed through the following
sections.

\begin{theorem} \label{thmMain1}
Consider a comparison constellation with lower tangency
$\{N^{n},S^{m},M_{w}^{m}\}$ on the interval
$[\,0, \infty[\,$. Assume further that the
functions $h(r)$ and $\lambda(r)$ are
{\em{balanced}} with respect to the warping
function $w(r)$ by the following inequality:
\begin{equation} \label{eqBalance1}
\mathcal{M}_p(r) :=  \left( m + p -
2\right)\eta_{w}(r) - m\,h(r) - (p-2)\lambda(r)
\geq  0.
\end{equation}
Let $\Lambda_{g,p}(r)$ denote the function
\begin{equation*}
\Lambda_{g,p}(r) =
w(r)\exp\left(-\int_{\rho}^{r}
\frac{\mathcal{M}_p(t)}{(p-1)g^{2}(t)}\,
dt\right).
\end{equation*}
Suppose finally that $p\ge 2$ and that
\begin{equation} \label{eqTransienceCond}
\int_{\rho}^{\infty} \Lambda_{g,p}(t)
\,  dt =  \infty.
\end{equation}
Then $S^{m}$ is $p$-parabolic.
\end{theorem}

\begin{theorem} \label{thmMain2}
Consider a comparison constellation with upper tangency
$\{N^{n},S^{m},M_{w}^{m}\}$ on the interval
$[\,0, \infty[\,$. Assume further that the
functions $h(r)$ and $\lambda(r)$ are
{\em{balanced}} with respect to the warping
function $w(r)$ by the following inequality:
\begin{equation} \label{eqBalance2}
\mathcal{M}_p(r) :=  \left( m + p -
2\right)\eta_{w}(r) - m\,h(r) - (p-2)\lambda(r)
\leq  0.
\end{equation}
Let $\Lambda_p(r)$ denote the function
\begin{equation*}
\Lambda_p(r) =
w(r)\exp\left(-\int_{\rho}^{r}
\frac{\mathcal{M}_p(t)}{(p-1)}\,
dt\right).
\end{equation*}
Suppose finally that $p\ge 2$ and that
\begin{equation} \label{eqTransienceCond2}
\int_{\rho}^{\infty} \Lambda_p(t)
\,  dt  = \infty.
\end{equation}
Then $S^{m}$ is $p$-parabolic.
\end{theorem}
\begin{remark}
It is easy to check that, when $p=2$, the lower bound $\lambda(r(x))$ for the $o$-radial component of the second fundamental form $\mathcal{B}(x)$ is obsolete when we consider a comparison constellation, (with upper or lower tangency), $\{N^{n},S^{m},M_{w}^{m}\}$ on the interval $[\,0, \infty[\,$. Moreover, the function $\mathcal{M}_2 (r)$ becomes in this case
$$\mathcal{M}_2(r)=m(\eta_w(r)- h(r)).$$
With this consideration at hand, we can find a version of Theorem \ref{thmMain1} for $p=2$ in the paper \cite{MP1}, where it is used a more restrictive balance condition which implies condition (\ref{eqBalance1}). On the other hand, we have, based on the same consideration, a version of Theorem \ref{thmMain2} for $p=2$ in the paper \cite{EP}, where we can find a direct proof and some consequences in connection with \cite{MP1}.
\end{remark}
\begin{corollary}\label{ThirdCor}
Consider a comparison constellation with upper tangency
$\{N^{n},S^{m},M_{w}^{m}\}$ on the interval $[\,0, \infty[\,$.
Assume that $\mathcal{M}_p(r) \leq 0$ for all $r>0$ and that the
warping function $w(r)$ is bounded from below by a positive
constant on $[r_0,\infty[$, for some $r_0>0$. Then, $S^m$ is
$p$-parabolic.
\end{corollary}
\begin{proof}
To obtain the result it suffices to apply Theorem \ref{thmMain2},
taking into account that, under the hypothesis, $\Lambda_p(r)\geq
w(r)$ for all $r>0$.
\end{proof}
%\begin{corollary}\label{FirstCor}
%Consider a comparison constellation with lower tangency
%$\{N^{n},S^{m},M_{w}^{m}\}$ on the interval $[\,0, \infty[\,$.
%Assume that, given $q \geq 2$,  $\mathcal{M}_q(r) \geq 0$ for all
%$r>0$ and $\eta_w (r)\geq \lambda(r)$ for all $r>0$.

%Assume furthermore that
%\begin{equation} \label{eqTransienceCond}
%\int_{\rho}^{\infty} \Lambda_{g,q}(t)
%\,  dt =  \infty.
%\end{equation}

%Then $S^{m}$ is $p$-parabolic, for all $p \geq q$
%\end{corollary}

\begin{corollary}\label{SecCor}
Consider a comparison constellation with upper tangency
$\{N^{n},S^{m},M_{w}^{m}\}$ on the interval $[\,0, \infty[\,$.
Assume further that,  given $q \geq 2$,  $\mathcal{M}_q(r) \leq 0$
for all $r>0$ and $h(r)\leq \eta_w (r)\leq \lambda(r)$ for all
$r>0$.

Assume furthermore that
\begin{equation} \label{eqTransienceCond2}
\int_{\rho}^{\infty} \Lambda_{q}(t)
\,  dt =  \infty.
\end{equation}

Then $S^{m}$ is $p$-parabolic, for all $p \geq q$.
\end{corollary}
\begin{proof}
It is straightforward to see that, if $p \geq q$,
\begin{equation}
\mathcal{M}_p(r) =\mathcal{M}_q(r)+(p-q)(\eta_w(r)
-\lambda(r))\leq \mathcal{M}_q(r)\leq 0,
\end{equation}
since $\eta_w(r)-h(r)<0$ for all $r>0$. Then, it is easy to check
that, under the hypothesis,
$$\frac{\mathcal{M}_p(r)}{p-1} \leq \frac{\mathcal{M}_q(r)}{q-1},$$

so $\int_{\rho}^{\infty} \Lambda_{p}(t) \, dt\geq
\int_{\rho}^{\infty} \Lambda_{q}(t) \,dt=\infty$.\\

Applying Theorem \ref{thmMain2}, the results follows.
\end{proof}

%%%%%%%%%%%%%%%%%%%%%%%%%%%%%%%%%%%%%%%%%%%%%%%%%%%
%%%%%%%%%%%%%%%%%%%%%%%%%%%%%%%%%%%%%%%%%%%%%%%%%%%%%%%%%%%%%%%%%%%%%%%%
%
%       SECTION 5 Proof of Theorem \ref{thmMain1}}
%
%%%%%%%%%%%%%%%%%%%%%%%%%%%%%%%%%%%%%%%%%%%%%%%%%%%%%%%%%%%%%%%%%%%%%%%%

%\section{Drifted $2$-capacity of model spaces} \label{secDrift}
\section{Proof of Theorem \ref{thmMain1}} \label{proofthmMain1}

We define on model spaces $M^m_w$, the modified Laplacian

\begin{displaymath}
\Lmod \psi(x) = \Delta^{M^m_w} \psi(x) + \psi'(r(x))\left(
\frac{\mathcal{M}(r(x))}{(p-1)\,g^{2}(r(x))} -
m\eta_{w}(r(x))\right),
\end{displaymath}
for smooth functions $\psi$ on $M^m_w$. If $\psi=\psi(r)$ only
depends on the radial distance r, then

\begin{equation}\label{operator}
\Lmod \psi(r) = \psi''(r) + \psi'(r)\left(
\frac{\mathcal{M}(r)}{(p-1)\,g^{2}(r)} -
\eta_{w}(r)\right).
\end{equation}

\noindent Consider now the following Dirichlet-Poisson problem
associated to $\Lmod$:

\begin{equation}\label{eqDir}
\begin{cases}
\Lmod \psi &= 0\,\,\,\text{on $A_{\rho, R}^{w}$},\\
\phantom{L }\psi &= 0\,\,\,\text{on $\partial B^w_\rho$}, \\
\phantom{L }\psi &= 1\,\,\,\text{on $\partial B^w_R$},
\end{cases}
\end{equation}
where $A_{\rho, R}^{w}$ is the annular domain in the model space defined as $A_{\rho, R}^{w}=B^w_R-B^w_\rho$.

The explicit solution to the Dirichlet problem
\eqref{eqDir} is given in the following
Proposition which is straightforward,

\begin{proposition}\label{propDirSol}
The solution to the Dirichlet problem (\ref{eqDir}) only depends on
$r$ and is given explicitly - via the function
$\Lambda_{g,p}(r)$ introduced in Theorem \ref{thmMain1}, by:
\begin{equation} \label{eqPsi}
\psi_{\rho,R}(r) =  \frac{\int_{\rho}^{r}
\Lambda_{g,p}(t)\,dt}{\int_{\rho}^{R}\Lambda_{g,p}(t)\,dt}.
\end{equation}
The corresponding 'drifted' $2$-capacity is
\begin{equation} \label{eqModelCap}
\begin{aligned}
\C_{\Lmod}(A_{\rho, R}^{w})
&=\int_{\partial D_{\rho}^{w}}\langle\nabla^{M}\psi_{\rho,R},\nu\rangle\,dA\\
&=\Vol(\partial D_{\rho}^{w})\Lambda_{g,p}(\rho)\left(\int_{\rho}^{R}
\Lambda_{g,p}(t)\,dt\right)^{-1}.
\end{aligned}
\end{equation}
\end{proposition}

It is easy to see, using equation (\ref{eqPsi}) and the balance condition (\ref{eqBalance1}) that
\begin{equation}
\begin{aligned}
\psi'_{\rho,R}(r)& \geq 0,\\
\psi''_{\rho,R}(r)& - \psi'_{\rho,R}(r)\eta_{w}(r)=-
\psi'_{\rho,R}(r)\frac{\mathcal{M}(r)}{(p-1)\,g^{2}(r)} \leq 0.
\end{aligned}
\end{equation}
 Now, we need the following result, which relates the $p$-Laplacian of a
radial function $f(r)$ with the operator  $\Lmod$.

\begin{lemma}\label{lemPLapComp}
Let $\{N^{n},S^{m},M_{w}^{m}\}$ be a comparison constellation with lower tangency  on $[0,R]$
for $R>0$. Let $f\circ r$ be a smooth real-valued function with $f' \geq 0$,
and suppose now that $f(r)$ satisfies the following
condition:
\begin{equation}\label{parenth}
f''(r) - f'(r)\eta_{w}(r) \leq  0.
\end{equation}
Then, for all $x \in S$,
\[
\pL^{S}f(r(x)) \leq (p-1)F^{p-2}(x)g^{2}(r(x))
\Lmod(f(r(x))) ,
\]
where $L$ is the second order differential operator defined
by equation (\ref{operator}) and $F$ is given by equation
\begin{equation}\label{eqabbrev}
F(x)=f'(r(x))\Vert \nabla^{S} r(x)\Vert.
\end{equation}
\end{lemma}
\begin{proof}
Computing as in \cite{HMP}, we have
\begin{equation}\label{eqLap0}
\begin{aligned}
&\pL^{S}f(r(x)) =
 F^{p-2}(x) \Bigl((p-2) \Bigl(
f''(r(x))\Vert \nabla^{S} r(x)\Vert^{2} \\&+
f'(r(x))\frac{\left\langle \nabla^{S}r(x), \nabla^{S}\Vert
\nabla^{S}r(x) \Vert\right \rangle}{\Vert \nabla^{S}r(x)
\Vert}\Bigr)+ \Delta^{S}f(r(x))\Bigr).
\end{aligned}
\end{equation}
This partial 'isolation' of the factor $(p-2)$ is
the reason behind the general assumption $p \geq
2$ in this work. Once we have equation (\ref{eqLap0}), we argue
as follows:

First, it is easy to see that
\begin{equation} \label{eqLambda}
\begin{aligned}
&\frac{\left \langle \nabla^{S}r(x), \nabla^{S}\Vert
\nabla^{S}r(x)\Vert\right \rangle }{\Vert \nabla^{S}r(x) \Vert}\,
\\ &\qquad \qquad \qquad =\Hess^{N}(r(x))\left(U_{r}, U_{r}\right)
+ \left\langle \nabla^{N}r(x),  \alpha_{x}\left(U_{r}, U_{r}
\right)\right \rangle.
\end{aligned}
\end{equation}
This quantity is bounded from above using Theorem \ref{thmGreW}
and the lower bound of $\mathcal{B}(x)$. On the other hand, since
the  $o$-radial mean convexity of $S$, $\mathcal{C}(x)$ is bounded
from below by the function $h(r(x))$, we obtain the following
estimate using Proposition  \ref{corLapComp}, (recall that
$f'(r)\geq 0$)
\begin{equation}
\label{eqLap1}
\Delta^{S}(f \circ r) \leq  \left(f''(r)
- f'(r)\eta_{w}(r) \right)
 \Vert \nabla^{S} r \Vert^{2}
+ m\,f'(r) \left(\eta_{w}(r) - h(r)\right).
\end{equation}

So, using the fact that $f(r)$ satisfies inequality (\ref{parenth}) and that  $\Vert \nabla^{S}(r)\Vert \geq
g(r)$, we have
\begin{equation}\label{eqLap2}
%\begin{aligned}
\pL^{S}(f(r(x))) \leq  (p-1)F^{p-2}(x)g^{2}(r) \Lmod(f(r)),
%\end{aligned}
\end{equation}
as claimed in the lemma.
\end{proof}

Now we transplant the model space
solutions $\psi_{\rho,R}(r)$ of equation
(\ref{eqDir}) into the extrinsic annulus
$A_{\rho,R}=D_{R}(o)\setminus \bar D_{\rho}(o)$ in
$S$ by defining
\[
\Psi_{\rho,R}\colon A_{\rho,R} \to \erre, \quad
\Psi_{\rho,R}(x)=\psi_{\rho,R}(r(x)).
\]

Since $\psi'_{\rho,R}(r) \geq 0$ and $\Lmod\psi_{\rho,R}=0$
in $A_{\rho, R}^{w}$, we obtain, applying Lemma \ref{lemPLapComp} to the function $\Psi_{\rho,R}$,

\begin{equation*}
\pLP\Psi_{\rho,R} \leq 0 \quad\text{in}\quad D_{R}(o)\setminus\bar
B_{\rho}(o),
\end{equation*}
that is to say,  $\Psi_{\rho,R}$ is a $p$-supersolution in
$D_{R}(o)\setminus\bar B_{\rho}(o)$. In fact, $\Psi_{\rho,R}$ is a
$p$-supersolution in the whole extrinsic ball $D_{R}(o)$ since
$\Psi_{\rho,R}(x)=0$ for $x\in S\cap\bar B_{\rho}(o)$; see
\cite{HMP}.

As $S$ is properly immersed, $D_{\rho}(o)$ and $D_{R}(o)$ are precompact and with regular boundary, so there exists a unique function $u\in C(\bar D_{R}(o))$
which is $p$-harmonic in $D_{R}(o)\setminus \bar D_{\rho}(o)$ such that
$u=0$ in $\bar D_{\rho}(o)$, $u=1$ in $\partial D_{R}(o)$, and that
\[
\pCap (\bar D_{\rho}(o),D_{R}(o))=\int_{D_{R}(o)}\Vert \nabla^{S} u\Vert^{p}\,d\mu.
\]
\noindent (see \cite{T1} and \cite[pp. 106-107]{HKM}).

Furthermore, let $\Psi_{\rho,R}$ be the transplanted $p$-supersolution
in $D_{R}(o)$. By the comparison principle, we have now
\[
u(x)\le\Psi_{\rho,R}(x)
\]
for all $x\in D_{R}(o)$. Hence, as $u(x)=\Psi_{\rho,R}(x)=0$ for
all $x\in \bar D_{\rho}(o)$, we obtain that
\begin{equation}\label{eqGradbound}
\Vert\nabla^{S}u(x)\Vert \le \Vert\nabla^{S}\Psi_{\rho,R}(x)\Vert
\end{equation}
for all $x\in \partial D_{\rho}(o)$.

With same arguments than in \cite{HMP}, but inverting all
inequalities, we obtain

\[
\pCap\bigl(\bar D_{\rho}(o),D_{R}(o)\bigr)
\le\left(\frac{\C_{\Lmod}(A_{\rho, R}^{w})} {\Vol(\partial
D_{\rho}^{w})}\right)^{p-1} \int_{\partial D_{\rho}}\Vert \nabla^S
r\Vert^{p-1}\, d\mathcal{H}^{m-1}.
\]
As, on the other hand,  $D_{\rho}(o)$ is precompact with a
smooth boundary thence,
\[
\int_{\partial D_{\rho}}\Vert \nabla^S
r\Vert^{p-1} \, d\mathcal{H}^{m-1} \, > 0 .
\]

So finally we have

\begin{equation}\label{intcond}
\begin{aligned}
&\pCap\bigl(\bar D_{\rho}(o), S^{m}\bigr)=\\ & \lim_{R\to \infty}
\pCap\bigl(\bar D_{\rho}(o),D_{R}(o)\bigr)\,\\& \leq
(\int_{\partial D_{\rho}}\Vert \nabla^S r\Vert^{p-1} \,
d\mathcal{H}^{m-1})\left(\lim_{R \to
\infty}\frac{\C_{\Lmod}(A_{\rho, R}^{w})} {\Vol(\partial
D_{\rho}^{w})}\right)^{p-1}=0,
\end{aligned}
\end{equation}
since  $\lim_{R \to \infty} \C_{\Lmod}(A_{\rho, R}^{w})=0$ by hypothesis
(\ref{eqTransienceCond}) and  equality (\ref{eqModelCap}) . Thus $\bar D_{\rho}(o)$ is a compact
subset with zero $p$-capacity in $S^{m}$, and $p$-parabolicity of
that submanifold follows.

%%%%%%%%%%%%%%%%%%%%%%%%%%%%%%%%%%%%%%%%%%%%%%%%%%%%%%%%%%%%%%%%%%%%%%%%
%
%       SECTION 6: Proof of Theorem \ref{thmMain2}
%
%%%%%%%%%%%%%%%%%%%%%%%%%%%%%%%%%%%%%%%%%%%%%%%%%%%%%%%%%%%%%%%%%%%%%%%%

\section{Proof of Theorem \ref{thmMain2}} \label{proofthmMain2}

We define the following modified Laplacian $\L$ on model spaces $M^m_w$, 
\begin{displaymath}
\L\phi(x) = \Delta^{M^m_w} \phi(x) + \phi'(r(x))\left(
\frac{\mathcal{M}(r(x))}{(p-1)} -
m\eta_{w}(r(x))\right),
\end{displaymath}
for smooth functions $\phi$ on $M^m_w$. As before, if $\phi=\phi(r)$ only
depends on the radial distance r, then

\begin{equation}\label{operator2}
\L\phi(r) = \phi''(r) + \phi'(r)\left(
\frac{\mathcal{M}(r)}{(p-1)} -
\eta_{w}(r)\right).
\end{equation}

\noindent Consider now the smooth radial solution $\phi_{\rho,R}(r)$ of the Dirichlet-Poisson problem associated to $\L$ and defined on the annulus $A_{\rho, R}^{w}=B^w_R-B^w_\rho$.

Now we transplant the model space
solutions $\phi_{\rho,R}(r)$ of this problem
into the extrinsic annulus
$A_{\rho,R}=D_{R}(o)\setminus \bar D_{\rho}(o)$ in
$S$ as in the proof  of Theorem \ref{thmMain1}, so we  have
\[
\Phi_{\rho,R}\colon A_{\rho,R} \to \erre, \quad
\Phi_{\rho,R}(x)=\phi_{\rho,R}(r(x)).
\]

Using the lower bounds of $\mathcal{B}(x)$, and $\mathcal{C}(x)$,
the fact that $\Phi_{\rho,R}'(r)\geq 0$ and applying Theorem
\ref{thmGreW}, and Proposition  \ref{corLapComp} we obtain, as we
did for any radial function $f(r)$ satisfying $f'(r)\geq 0$ in the
proof of Lemma \ref{lemPLapComp}:
\begin{equation}
\begin{aligned}
\Delta^{S}\Phi_{\rho,R} &\leq  \left(\Phi_{\rho,R}''(r)
- \Phi_{\rho,R}'(r)\eta_{w}(r) \right)
 \Vert \nabla^{S} r \Vert^{2}
\\ &+ m\,\Phi_{\rho,R}'(r) \left(\eta_{w}(r) - h(r)\right).
\end{aligned}
\end{equation}

Hence, as
$$\Phi_{\rho,R}''(r)
- \Phi_{\rho,R}'(r)\eta_{w}(r) \,\, \geq \,\, 0 \,\,\forall r>0$$
because the balance condition  (\ref{eqBalance2}) and $\Vert
\nabla^S r\Vert \leq 1$, we obtain
\begin{equation}
\pL^{S}(\Phi_{\rho,R}(r(x))) \leq  (p-1)F^{p-2}(x)g^{2}(r) \Lmod(\Phi_{\rho,R}(r))=0,
\end{equation}
where $F=\Phi_{\rho,R}'(r(x))\Vert\nabla^S r(x)\Vert$. The rest of
the proof follows in the same way than in Section
\ref{proofthmMain1}.

%%%%%%%%%%%%%%%%%%%%%%%%%%%%%%%%%%%%

%\bibliographystyle{amsalpha}

%\end{document}

\enddocument